\input ./style/arxiv-general.cfg
\documentclass[aos,MSNbibl,nameyear,dvips]{arximspdf}
\makeatletter
   \@ifpackageloaded{graphicx}{}{\usepackage{graphicx}}
\makeatother
\usepackage{dcolumn}

\doi{10.1214/15-AOS1370}
\volume{44}
\issue{1}
\pubyear{2016}
\firstpage{401}
\lastpage{424}
\docsubty{FLA}

\makeatletter
\newcolumntype{d}[1]{D{.}{.}{#1}}

\newcommand{\rrvert}{\vert}
\newcommand{\rrVert}{\Vert}
\newcommand{\llvert}{\vert}
\newcommand{\llVert}{\Vert}

\newcommand{\eqref}[1]{(\ref{#1})}
\newtheorem{theorem}{Theorem}[section]
\newtheorem{lemma}[theorem]{Lemma}
\newtheorem{corollary}[theorem]{Corollary}
\newproclaim{definition}{Definition}
\newproclaim{remark}{Remark}
\makeatother

\begin{document}
\begin{frontmatter}

\title{A goodness-of-fit test for stochastic block models}
\runtitle{A goodness-of-fit test for stochastic block models}

\begin{aug}
\author[A]{\fnms{Jing}~\snm{Lei}\corref{}\ead[label=e1]{jinglei@andrew.cmu.edu}\ead[label=u1,url]{http://www.stat.cmu.edu/\textasciitilde jinglei}}
\runauthor{J. Lei}
\affiliation{Carnegie Mellon University}
\address[A]{Department of Statistics\\
Carnegie Mellon University\\
Pittsburgh, Pennsylvania 15213\\
USA\\
\printead{e1}}
\end{aug}

%
\received{\smonth{12} \syear{2014}}
%
\revised{\smonth{8} \syear{2015}}

%
\begin{abstract}
The stochastic block model is a popular tool for
studying community structures in network data.
We develop a goodness-of-fit test for the stochastic block model.
The test statistic is based on the largest singular value
of a residual matrix obtained by subtracting the
estimated block mean effect from the adjacency matrix.
Asymptotic null distribution is obtained using recent
advances in random matrix theory. The test is proved to
have full power against alternative models with finer structures.
These results naturally lead to a consistent sequential testing estimate
of the number of communities.
\end{abstract}

%
\end{frontmatter}

\section{Introduction}
\label{sec:introduction}
Large-scale network data with community structures have been the focus
of much research efforts in the past decade
[see, e.g., \citet{NewmanG04,Newman06}].
In the statistics and machine learning literature,
the stochastic block model \citep{Holland83} is a very
popular model for community structures in network data.
In a stochastic block model, the observed network is
often recorded in the form
of an $n\times n$ adjacency matrix $A$, representing
the presence/absence of pairwise interactions
among $n$ individuals in a population of interest.
The model assumes that (i) the individuals are partitioned into
$K$ disjoint communities, and (ii) given the memberships,
the upper diagonal entries of $A$ are independent Bernoulli
random variables, where the parameter $E(A_{ij})$ depends only on the
memberships
of nodes $i$ and $j$.
Such a
model naturally captures the community structures commonly observed
in complex networks, and has close connection to nonparametric
exchangeable random graphs \citep{BickelC09}.
The stochastic block model can be made more realistic by
incorporating additional parameters to better approximate real world network
data. For example, \citet{KarrerN11} incorporated individual node
activeness into the stochastic block model to allow for arbitrary
degree distributions.
In the mixed membership model \citep{Airoldi08}, each individual
may belong to more than one community.

In this paper, we develop a goodness-of-fit test for stochastic
block models.
Given an adjacency matrix $A$ and a positive integer $K_0$,
we test whether $A$ can be adequately fitted by a stochastic block
model with $K_0$
communities.
Our test statistic is the largest singular value of a residual matrix
obtained by removing the estimated block mean effect from the observed
adjacency matrix.
Intuitively, if $A$ is generated by a stochastic block model and the
block mean
effect is estimated appropriately, the residual matrix will approximate
a generalized Wigner matrix:
a symmetric random matrix with independent mean zero upper diagonal entries.
Our first contribution is the asymptotic null distribution of the test
statistic (Theorem~\ref{thm:asymp-dist}).
The proof uses some recent advances in random matrix theory, such as
the local
semicircle law of generalized Wigner matrices and its consequences
\citep{Rigidity,General,Isotropic}.
Our second contribution is asymptotic power guarantee of the test against
models with finer structures (Theorems~\ref{thm:consistency} and~\ref{thm:dcbm}).
In particular, we establish the growth rate
of the test statistic under alternatives that correspond to stochastic
block models
with more communities or with individual node degree variations. It is
of particular interest
to consider alternative stochastic block models with more communities
because any exchangeable random graph
can be approximated by a stochastic block model \citep{BickelC09}. In
our simulation study, we observe that the proposed test is powerful
against not only stochastic block
models with more communities, but also other network models with finer
structures such as the degree corrected
block model and the mixed membership block model.

A related test statistic using the largest eigenvalue of the centered
and scaled adjacency
matrix has been studied in \citet{BickelS13test}.
They derive asymptotic null distribution for Erd\H{o}s--R\'{e}nyi models,
which corresponds to a stochastic block model with one community.
We generalize their argument to prove the asymptotic null distribution
result for stochastic block models with more
than one community. The key step is to bound the fluctuation in
the leading eigenvalue of a random matrix under perturbation of a
block-wise constant
noise matrix.
Moreover, their asymptotic power analysis requires the alternative
model to be diagonal dominant.
Our test statistic uses the largest singular value of the residual
matrix, so we are able to capture signals affecting either the largest
or the smallest
eigenvalues, and our asymptotic power guarantee holds for a much wider
class of
alternative models.


Our goodness-of-fit test can also serve as
a main building block to estimate the number of communities.
As a key inference problem in stochastic block models and its variants,
the community recovery problem concerns estimating the hidden
communities from
a single observed adjacency matrix
[see \citet{McSherry01,BickelC09,Decelle11,ZhaoLZ12,Jin12,Fishkind13,LeiR14,ChenSX12,ChaudhuriCT12,Krzakala13,Massoulie13,MosselNS13,AbbeBH14,Anandkumar14}, e.g.].
A common assumption made in all these methods is that $K$, the total
number of communities, is known.
Therefore, estimating the number of communities in a stochastic block
model is
of great practical and theoretical importance.
Some methods have been proposed to estimate the number of communities
in stochastic block models
\citep{ZhaoLZ11PNAS,BickelS13test,ChenL14,SaldanaYF14}, but without
consistency guarantee.

To estimate the number of communities, we consider hypothesis test
%
\begin{equation}
\label{eq:hypothesis} H_{0,K_0}: K=K_0,\quad\mbox{against}\quad
H_{\mathrm{a},K_0}: K>K_0
\end{equation}
sequentially
for each $K_0\ge1$ until the null hypothesis is not rejected.
We prove the consistency of this sequential testing estimator in
Corollary~\ref{cor:consistency} of
Section~\ref{sec:theory}. Throughout this paper, we use $K$ to denote
the true number of communities
in a stochastic block model and use $K_0$ to denote a hypothetical
number of communities.

Recently, \citet{Chatterjee15} studied a general method for matrix denoising
using singular value thresholding, which covers the stochastic block
model as a special case.
Following the ideas developed there, one may use the number of
significant singular values, for example, those greater than $\sqrt
{n}$, of the adjacency matrix as an estimate of $K$. But this method only
works when the community-wise connectivity matrix has full rank.
Empirically, we also find that the sequential
testing estimator developed in this paper performs better than singular
value thresholding for sparse networks.

\subsection*{Glossary} For a square matrix $M$, $\operatorname{diag}(M)$
denotes the diagonal matrix induced by $M$. For any $n\times n$
symmetric matrix $M$, $\lambda_j(M)$ denotes
its $j$th largest eigenvalue value, ordered as $\lambda_1(M)\ge\lambda
_2(M)\ge\cdots\ge\lambda_n(M)$, and $\sigma_1(M)$ is the largest
singular value.
Denote $\mathcal B_K$ the set of all $K\times K$
symmetric matrices with entries in $(0,1)$ and all rows being distinct.

\section{Stochastic block models and a goodness-of-fit test}
\label{sec:sbm-seq}
A stochastic block model on $n$ nodes with $K$ communities is
parameterized by
a membership vector $g\in\{1,\ldots,K\}^n$ and a symmetric community-wise
edge probability matrix $B\in[0,1]^{K\times K}$.
The observed adjacency matrix $A$ is a symmetric binary matrix with diagonal
entries being 0.
Given $(g,B)$, the
probability mass function for the adjacency matrix $A$ is
%
\begin{equation}
\label{eq:sbm} P_{g,B}(A)=\prod_{1\le i< j\le n}
B_{g_i g_j}^{A_{ij}} (1-B_{g_i
g_j})^{(1-A_{ij})}.
\end{equation}
In other words, given $(g,B)$, the edges are independent Bernoulli random
variables with parameters determined by the node memberships.

To avoid triviality, we say that a stochastic block model parameterized
by $(g,B)$
has $K$ communities if (i) $g$ contains all $K$ distinct values
in $\{1,\ldots,K\}$, and (ii)~any two rows of $B$ are distinct.
A stochastic block model is identifiable up to
a label permutation on $g$ and a corresponding row/column permutation
on $B$.

Given an observed adjacency matrix $A$,
and a positive integer $K_0$,
we would like to know if $A$ can be well fitted by a stochastic block
model with $K_0$
communities.
If we assume that $A$ is generated by a stochastic block model with $K$
communities,
this leads to a goodness-of-fit test for stochastic block models
with a composite null hypothesis
%
\begin{equation}
\label{eq:hypothesis0} H_{0,K_0}: K=K_0.
\end{equation}

To derive a goodness-of-fit test for stochastic block models,
a natural idea is to estimate the model parameters and remove the signal
from the observed adjacency matrix, and test whether the residual matrix
looks like a noise matrix.
To this end,
consider the $n\times n$ matrix $P$ given by
\[
P_{ij}=B_{g_i g_j}, %
\]
so that $E(A)=P-\operatorname{diag}(P)$.
Let $\tilde A^*$ be
%
\begin{equation}
\label{eq:A*} \tilde A^*_{ij}= \frac{A_{ij}-P_{ij}}{\sqrt{(n-1)P_{ij}(1-P_{ij})}},\qquad i\neq j
\quad\mbox{and}\quad \tilde A^*_{ii}=0, \forall i.
\end{equation}
Then $\tilde A^*$ is a generalized Wigner matrix, satisfying
$E(\tilde A^*_{ij})=0$ for all $(i,j)$ and $\sum_{j}\operatorname{var}(\tilde
A^*_{ij})=1$ for all $i$.
The asymptotic distribution
of the extreme eigenvalues of $\tilde A^*$ has been well studied in
random matrix theory.
In particular, combining results in \citet{Rigidity} and
\citet{LeeY14} we have
%
\begin{equation}
\label{eq:A*TW} 
n^{2/3} \bigl[\lambda_1\bigl(
\tilde A^*\bigr)-2 \bigr]\rightsquigarrow TW_1\quad\mbox{and}\quad
n^{2/3} \bigl[-\lambda_n\bigl(\tilde A^*\bigr)-2 \bigr]
\rightsquigarrow TW_1, 
\end{equation}
where $TW_1$ denotes the Tracy--Widom distribution with index 1
and ``$\rightsquigarrow$'' denotes convergence in distribution. We remark
that \eqref{eq:A*TW} cannot be obtained using results for standard
Wigner matrices
as the diagonal entries of $\tilde A^*$ are fixed to be 0. We formally
state and prove this
result as Lemma~\ref{lem:TW-A*} in Section~\ref{sec:RMT}.


The matrix $\tilde A^*$ involves unknown model parameters and cannot be
used as a
test statistic. Now we describe a natural estimate of $\tilde A^*$ by
plugging in
an estimated stochastic block model.

Let $\hat g$ be an estimated community membership vector with target
number of
communities being $K_0$.
Define $\hat{\mathcal N}_{k}=\{i:1\le i\le n, \hat g_i=k\}$, and $\hat
n_k=\llvert  \hat{\mathcal N}_k\rrvert  $ for all $1\le k\le K_0$.
We consider the plug-in estimator of $B$:
%
\begin{equation}
\label{eq:plug-in-B} \hat B_{kl} = \cases{ \displaystyle\frac{\sum_{i\in\hat{\mathcal N}_k,j\in\hat{\mathcal N}_l}
A_{ij}}{\hat n_k \hat n_l}, &
\quad$k\neq l$,
\cr
\displaystyle\frac{\sum_{i,j\in\hat{\mathcal N}_k,i<j}A_{ij}}{\hat
n_k(\hat n_k-1)/2},&\quad$k=l $.}
\end{equation}
The estimates $(\hat g,\hat B)$ leads to the empirically centered and re-scaled
adjacency matrix $\tilde A$:
%
\begin{equation}
\label{eq:tilde-A} \tilde A_{ij} =\cases{ \displaystyle \frac{A_{ij}-\hat P_{ij}}{\sqrt{(n-1)\hat P_{ij}(1-\hat
P_{ij})}},&
\quad$i\neq j$,
\vspace*{4pt}\cr
\displaystyle0,&\quad$i=j$,}
\end{equation}
where
%
\begin{equation}
\label{eq:est-P} \hat P_{ij}=\hat B_{\hat g_i\hat g_j}.
\end{equation}
It is natural to conjecture that under the null hypothesis $K=K_0$ and
when the estimates $(\hat g,\hat B)$ are accurate enough,
the convergence in \eqref{eq:A*TW} will carry over to
the corresponding eigenvalues of $\tilde A$. Therefore,
we can use the largest singular value of $\tilde A$, after centering
and scaling, as our test statistic:
%
\begin{equation}
\label{eq:test-stat} T_{n,K_0}=n^{2/3} \bigl[\sigma_1(
\tilde A)-2 \bigr].
\end{equation}
The corresponding level $\alpha$ rejection rule for testing problem
\eqref{eq:hypothesis0} is
%
\begin{equation}
\label{eq:rule} \mbox{Reject}\qquad H_{0,K_0},\qquad\mbox{if }
T_{n,K_0}\ge t(\alpha/2),
\end{equation}
where $t(\alpha/2)$ is the $\alpha/2$ upper quantile of the $TW_1$ distribution
for $\alpha\in(0,1)$. We use $t(\alpha/2)$ instead of $t(\alpha)$ for
Bonferroni correction because
\[
\sigma_1(\tilde A) = \max\bigl(\lambda_1(\tilde A),-
\lambda_n(\tilde A)\bigr), %
\]
and hence
\[
T_{n,K_0}=\max \bigl[ n^{2/3}\bigl(\lambda_1(\tilde
A)-2\bigr), n^{2/3}\bigl(-\lambda_n(\tilde A)-2\bigr) \bigr].
\]

A similar result
concerning the largest eigenvalue of $\tilde A$ in the simple case
\mbox{$K_0=1$} has
been obtained in \citet{BickelS13test}. In Section~\ref{sec:theory}
below, we
formally state and prove
the validity of our test statistic $T_{n,K_0}$ in Theorem~\ref
{thm:asymp-dist} by establishing
the asymptotic null distributions of
both the largest and smallest eigenvalues and for general
values of $K_0$.

The use of $\sigma_1(\tilde A)$ instead of $\lambda_1(\tilde A)$ as our
test statistic
in \eqref{eq:test-stat} is crucial for power guarantee. Under some
alternative hypotheses,
the signal may be carried solely by $\lambda_n(\tilde A)$. For example,
consider a model with two equal-sized communities and
$B_{11}=B_{22}=1/4$, $B_{12}=B_{21}=1/2$.
Suppose we would like to test $H_0:  K=  K_0=1$.
In this case, $A-\hat P$ has block-wise mean
\[
\pmatrix{ -1/8 & 1/8
\cr
1/8 & -1/8},
\]
which has no positive eigenvalues. Therefore,
the test using only $\lambda_1(\tilde A)$ has no power.

Given the rejection rule \eqref{eq:rule} for testing problem \eqref
{eq:hypothesis0},
we have the following sequential testing estimator of $K$:
%
\begin{equation}
\label{eq:K-est} \hat K = \inf\{K_0\ge1: T_{n,K_0}<
t_n\}.
\end{equation}
In other words, we perform the goodness-of-fit test for
$K_0=1,2,\ldots,$ until failing to reject $H_{0,K_0}$. We prove
consistency of $\hat K$ for appropriate choices of $t_n$
in Corollary~\ref{cor:consistency} below, as a consequence of
(i)
a large deviation inequality of the extreme eigenvalues of $\tilde A$
under the
null hypothesis $K=K_0$,
and (ii) the growth rate of $T_{n,K_0}$ under the alternative
hypothesis $K>K_0$.

\section{Asymptotic null distribution and power}
\label{sec:theory}
The asymptotic distribution of the test statistic
$T_{n,K_0}$ under the null hypothesis
depends on the accuracy of the estimated community membership
$\hat g$.
In order to consider the asymptotic behavior of community recovery,
we consider\vspace*{1pt} a sequence of stochastic block models
$\{(g^{(n)},B^{(n)}):n\ge1\}$ where $g^{(n)}\in\{1,\ldots,K^{(n)}\}^n$
for each $n$,
and $B^{(n)}\in\mathcal B_{K^{(n)}}$. Here, the number of communities
$K=K^{(n)}$ and
the community-wise edge probability matrix $B=B^{(n)}$ are allowed to
change with $n$.

We will focus on
relatively balanced communities.
\begin{enumerate}[(A1)]
\item[(A1)] There exists $c_0>0$ such that
$\min_{1\le k\le K^{(n)}}\llvert  \{i:g_i^{(n)}=k\}\rrvert   \ge c_0 n/K^{(n)}$
for all~$n$.
\end{enumerate}
Assumption (A1) assumes each community has size at least proportional
to $n/K^{(n)}$.
For example, it is satisfied almost surely if the membership vector
$g^{(n)}$ is
generated from a multinomial distribution with $n$ trials and
probability $\pi=(\pi_1,\ldots,\pi_{K^{(n)}})$ such that $\min_{1\le
k\le K}\pi_k>c_0/ K^{(n)}$
and $K^{(n)}$ grows slowly.



%
%

\begin{definition*}[(Consistency of community recovery)]
For\vspace*{1pt} a sequence of stochastic block models $\{(g^{(n)},B^{(n)}):n\ge1\}
$ with $K^{(n)}$ communities and
$B^{(n)}\in\mathcal B_{K^{(n)}}$,
we say a community membership estimator $\hat g=\hat g(A,K^{(n)})$ is consistent
if
\[
P_{A\sim(g^{(n)},B^{(n)})}\bigl(\hat g = g^{(n)}\bigr)\rightarrow1.
\]
\end{definition*}

\begin{remark*}
The notion ``$\hat g=g$'' shall be interpreted as being equal up to a label
permutation. Such a label permutation does not affect our methodological
and theoretical development so we will assume that the label permutation
is identity for simplicity.
The definition of consistent community recovery can be satisfied
by several methods. For example, in the case of fixed finite
$K^{(n)}=K$ and \mbox{$B^{(n)}=B$}, the profile likelihood method \citep{BickelC09}
is consistent for all $(g^{(n)}:n\ge1)$ satisfying (A1) and all
$B\in\mathcal B_K$; the spectral clustering method
can be made consistent, with slight modification, for all
$(g^{(n)}:n\ge1)$ satisfying (A1) and $B\in\mathcal B_K$
with full rank \citep{McSherry01,Vu14,LeiZ14}.
In the case of slowly growing $K^{(n)}$, consistent community
recovery can be achieved in some special cases such
as the planted partition model \citep{ChaudhuriCT12,AminiL14}.
\end{remark*}

\subsection{The asymptotic null distribution}
%
\begin{theorem}[(Asymptotic null distribution)]\label{thm:asymp-dist}
Let $A$ be an adjacency matrix generated from stochastic block model
$(g^{(n)},B^{(n)})$, where $B^{(n)}\in\mathcal B_{K^{(n)}}$ and
$(g^{(n)}:n\ge1)$ satisfies
condition \textup{(A1)}. Let $\tilde A$ be given as in \eqref{eq:tilde-A} using
a consistent community estimate $\hat g$ and corresponding plug-in estimate
of $B^{(n)}$ as in \eqref{eq:plug-in-B}.
Assume
in addition that $K^{(n)}=O(n^{1/6-\tau})$ for some $\tau>0$, and the
entries of $B^{(n)}$
are uniformly bounded away from 0 and 1.
The following holds under the null hypothesis $K^{(n)}=K_0$:
%
\begin{equation}
\label{eq:TW-A-tilde} n^{2/3}\bigl(\lambda_1(\tilde A)-2\bigr)
\rightsquigarrow TW_1,\qquad n^{2/3}\bigl(-
\lambda_n(\tilde A)-2\bigr)\rightsquigarrow TW_1.
\end{equation}
\end{theorem}

Theorem~\ref{thm:asymp-dist} is proved in Section~\ref{sec:proof-null}.
The main challenge
is that, assuming $\hat g=g$, the entry-wise estimation error in $\hat
B$ is of order
$K^{(n)}/n$. The simple upper bound of $\tilde A-\tilde A^*$ in
Frobenius norm is
of order $K^{(n)} n^{-1/2}$ which exceeds the $n^{-2/3}$ scaling
required in \eqref{eq:TW-A-tilde}.
Our proof establishes \eqref{eq:TW-A-tilde} using a more delicate analysis
that exploits the block-wise constant structure in $\tilde A-\tilde
A^*$, combined
with random matrix theory results which ensure that (i) the leading
eigenvectors of $\tilde A^*$
are delocalized, in the sense that the chance these eigenvectors being
close to any fixed vector is small, and (ii) the number of large
eigenvalues of $\tilde A^*$ in an interval
of length $K^{(n)}/\sqrt{n}$ can be accurately
approximated. This result is a nontrivial generalization of Theorem 2.1
in \citet{BickelS13test}.

An immediate consequence of Theorem~\ref{thm:asymp-dist} is an
asymptotic type I error bound
for the rejection rule \eqref{eq:rule}:
\begin{eqnarray*}
&&P \bigl[T_{n,K_0}\ge t(\alpha/2) \bigr]
\\
&&\qquad \le P \bigl[n^{2/3}\bigl(\lambda _1(\tilde A)-2\bigr)\ge
t(\alpha/2) \bigr] +P \bigl[n^{2/3}\bigl(-\lambda_n(\tilde
A)-2\bigr)\ge t(\alpha/2) \bigr]
\\
&&\qquad =\alpha/2+o(1)+\alpha/2+o(1)=\alpha+o(1).
\end{eqnarray*}
Formally, we have the following corollary.

\begin{corollary}[(Asymptotic type I error control)]\label{cor:level}
Under the assumptions of Theorem~\ref{thm:asymp-dist}, the rejection
rule in \eqref{eq:rule} has asymptotic level $\alpha$.
\end{corollary}

\subsection{Asymptotic power against $K>K_0$ and consistent estimation of $K$}
Now we consider the power of the test against finer stochastic block models.
The following theorem provides a lower bound of the growth rate of the
test statistic $T_{n,K_0}$ under the alternative model $K^{(n)}>K_0$.

\begin{theorem}[(Asymptotic power guarantee)]\label{thm:consistency}
Let $A$ be an adjacency matrix generated from stochastic block model
$(g^{(n)},B^{(n)})$ with
$B^{(n)}\in\mathcal B_{K^{(n)}}$ and $(g^{(n)}:n\ge1)$ satisfying
condition \textup{(A1)}.
Let $\delta_n$ be the smallest $\ell_\infty$ distance
among all pairs of distinct rows of $B^{(n)}$.
For any $K_0<K^{(n)}$ and any community estimator $\hat g$, we have
\[
\sigma_1(\tilde A)\ge\tfrac{1}{2}\delta_n
c_0 \bigl[K^{(n)}\bigr]^{-2} n^{1/2}+O_P(1).
\]
\end{theorem}

Theorem~\ref{thm:consistency} is powerful in that it puts no structural
condition on
the connectivity matrix $B^{(n)}$, nor does it make any assumption
about the particular
method used to estimate the membership.
Theorem~\ref{thm:consistency} is proved in Section~\ref
{sec:proof-power}. The main idea is that
if the nodes are partitioned into less than $K^{(n)}$ groups, the
corresponding block
partition of the expected adjacency matrix cannot be block-wise
constant, and hence
it is impossible to remove the mean effect by subtracting a constant
from each estimated
block submatrix of $A$.

When $B^{(n)}=B$ and $K^{(n)}=K$ are fixed and do not change with $n$,
the community separation parameter $\delta_n$ is constant
and Theorem~\ref{thm:consistency} gives a growth rate at least
$n^{1/2}$. When $K^{(n)}$ is allowed to grow with $n$,
consistent community recovery can be achieved for the planted partition
model where $B^{(n)}_{kk}=p$ and $B^{(n)}_{kk'}=q$ ($k\neq k'$) for
some $0\le q< p\le1$.
If $p$ and $q$ are constants independent of $n$, then $\delta_n=p-q$ is
also a constant. Therefore,\vspace*{1pt} in this case Theorem~\ref{thm:consistency}
says that the growth rate of $T_{n,K_0}$ is at least $[K^{(n)}]^{-2} n^{1/2}$.

The asymptotic null distribution and growth rate under alternative $K^{(n)}>K_0$
suggest that the null and alternative hypotheses are well separated. Therefore,
if in the sequential testing estimator \eqref{eq:K-est} we choose the
rejection threshold $t_n$ to increase with the network size $n$,
we shall expect to have a consistent estimate of~$K^{(n)}$.

\begin{theorem}[(Consistency of estimating $K$)]\label{cor:consistency}
Under the assumptions of Theorem~\ref{thm:asymp-dist} and Theorem~\ref
{thm:consistency}, assume in addition that
$\lim\inf_{n\rightarrow\infty}\delta_n > 0$. Let $\hat K$ be the
sequential testing estimator given in \eqref{eq:K-est}
with threshold $t_n$ satisfying $t_n\asymp n^{\varepsilon}$ for some
$\varepsilon\in(0,5/6)$, then
\[
P\bigl(\hat K=K^{(n)}\bigr)\rightarrow1. %
\]
\end{theorem}

Corollary~\ref{cor:consistency} is proved in Section~\ref{sec:proof-power}.
We note that the asymptotic null distribution given in Theorem~\ref
{thm:asymp-dist} cannot
be directly used to bound the probability of $P(T_{n,K^{(n)}}\ge t_n)$
because $t_n$ changes
with $n$. Instead, we need to use a tail probability bound on the
largest singular
value of $\tilde A$ (Lemma~\ref{lem:deviation}).
The condition that $\delta_n$ is bounded away from zero
is satisfied both when $B^{(n)}$ is fixed or when $B^{(n)}$ is given
by a planted partition model with constant diagonal and off-diagonal
edge probabilities.
This condition can be relaxed to requiring $\delta_n$ to decay no
faster than $n^{-1/6}$
and having $t_n\ll n^{5/6}\delta_n$.

\subsection{Asymptotic power against degree corrected block models}
The good\-ness-of-fit test \eqref{eq:rule} is also powerful against
certain degree corrected block models.
A~degree corrected block model is parameterized by a triplet $(g,B,\psi
)$, where
$\psi\in[0,1]^n$ is the node activeness parameter and the edge probability
between nodes $i$ and $j$ is $\psi_i\psi_j B_{g_i g_j}$.
The probability mass function of the observed adjacency matrix is
\[
P_{g,B,\psi}(A)=\prod_{1\le i<j\le n} (\psi_i
\psi_j B_{g_i g_j})^{A_{ij}}(1-\psi_i
\psi_j B_{g_i
g_j})^{1-A_{ij}}.
\]
%
Let $\phi_k$ be the subvector of $\psi$ corresponding to the entries in
community $k$,
and $\tilde\phi_k = \phi_k/\llVert  \phi_k\rrVert  $.

The condition we will need on $\psi$ is that there exists a community
whose node activeness parameters cannot be approximated by block-wise
constant vectors.
Formally, for any vector $v$ and positive integer $L$ let
\[
\mathcal E(v,L)=\min\bigl\{\llVert v-u\rrVert _2^2:
\mbox{the entries of } u \mbox{ take at most } L \mbox{ distinct values}\bigr\}.
\]

%
\begin{theorem}\label{thm:dcbm}
Let $A$ be generated by a degree corrected block model $(g,B,\psi)$ on
$n$ nodes and $K$ communities. If there exists $1\le k^*\le K$ such that
$\mathcal E(\tilde\phi_{k^*},K_0)>0$, then for any estimator $(\hat
g,\hat B)$ of a $K_0$-stochastic block model, we have
\[
\llVert \tilde A\rrVert \ge\frac{\mathcal E(\tilde\phi_{k^*},K_0)}{2K_0^{3/2}}\llVert B_{k^*,\cdot}\rrVert
_\infty\kappa_n n^{1/2}+O_P(1),
\]
where $\kappa_n=\min_{1\le k\le K}\llVert  \phi_k\rrVert  ^2 / n $ and $\llVert  B_{k^*,\cdot
}\rrVert  _\infty=\max_{1\le k\le K}B_{k^*,k}$.
\end{theorem}

Theorem~\ref{thm:dcbm} is proved in Section~\ref{sec:proof-dcbm}.
The quantity $\mathcal E(\tilde\psi_{k^*},K_0)$ reflects the idea that
there exists at least one community
whose node activeness cannot be approximated by a simple $K_0$-block structure.
Recall that for each $k$, $\tilde\phi_k$ is a vector with unit $\ell
_2$ norm. If the entries of $\phi_k$ are sampled from a compact
interval with strictly positive density,
then $\mathcal E(\tilde\phi_k, K_0)\asymp K_0^{-2}$ when $K_0$ is
small compared to the length of $\phi_k$.
When $K_0$ increases, $\mathcal E(\tilde\phi_k,K_0)$ decreases for all
$k$ and the test will
be less powerful. This agrees with the fact that any degree corrected
block model can
be approximated by regular stochastic block models with a large number
of communities.

Consider an opposite case, where $A$ is generated by a degree corrected
block model with one community and degree parameter vector $\psi$
containing only $K_0$ distinct values. Here, the model can also be
viewed as a regular stochastic block model with $K_0$ communities. Then
$\mathcal E(\tilde\psi_{k^*},K_0)=0$ and the test will not tend to
reject the null hypothesis, provided that a consistent community
recovery method is used.

The quantity $\kappa_n$ acts as a lower bound on the overall node
activeness. A larger value
of $\kappa_n$ usually leads to better power because
there are more observed edges for inference.
Under the balanced community assumption (A1), $\kappa_n\asymp K^{-1}$
if the entries of $\psi$ are uniformly bounded away from zero, or are
sampled from a common distribution independent of $n$.

Applying Theorem~\ref{thm:dcbm} in the simple special case where $B\in
\mathcal B_K$ (and hence $K$) is fixed and $\psi_i$'s are sampled from
a compact interval with strictly positive density, under Assumption
(A1) we have, for any given $K_0$,
\[
\llVert \tilde A\rrVert \ge C n^{1/2}+O_P(1). %
\]
Therefore, with probability tending to one, the test will reject the
null hypothesis that
$A$ is generated from a regular stochastic block model with $K_0$ blocks.
If $K$ grows with $n$ and the entries of $B$ scale at rate $\rho_n$,
the test is still powerful as long as
$n^{1/2}\rho_n/(K_0^{7/2}K)\rightarrow\infty$.

\section{Numerical experiments}
\label{sec:experiment}
Now we illustrate the performance of the proposed test and the
estimator of $K$ in various
simulations. In our simulation, we use simple spectral clustering for
community recovery.
Given an adjacency matrix $A$ and a hypothetical number of communities
$K_0$, this algorithm estimates the community membership by applying
$k$-means clustering
to the rows of the matrix formed by the $K_0$ leading singular vectors
of $A$.

\subsection{Simulation 1: The null distribution and bootstrap correction}
In the first simulation, we consider the finite sample null
distribution of the
scaled and centered extreme eigenvalues of $\tilde A$ and empirically verify
Theorem~\ref{thm:asymp-dist} for a simple stochastic block model.
Following the observation in \citet{BickelS13test}, the
speed of convergence to the limit distribution may be slow. A practical
solution to this issue
using a fused bootstrap correction has been
proposed in \citet{BickelS13test} for the special case of $K_0=1$. Here, we
extend this idea to the more general case considered in this paper.

For a given adjacency matrix $A$ on $n$ nodes and null hypothesis
$K=K_0$, the goodness-of-fit test statistic with fused bootstrap
correction is given as follows:
\begin{enumerate}[4.]
\item[1.] Let $\hat g$ be an estimated community membership vector with $K_0$ communities,
and $(\hat B,\hat P)$ be the corresponding estimates in \eqref
{eq:plug-in-B} and \eqref{eq:est-P}.\vspace*{1pt}

\item[2.] Calculate $\tilde A$ as in \eqref{eq:tilde-A} and its extreme
eigenvalues $\lambda_1(\tilde A)$, $\lambda_n(\tilde A)$.
\item[3.] For $m=1,\ldots,M$:
\begin{longlist}[(a)]
\item[(a)] Let $A^{(m)}$ be an adjacency matrix independently generated from
stochastic block model $(\hat g,
\hat B)$.
\item[(b)] Let $\tilde A^{(m)}=(\tilde A^{(m)}_{ij})_{i,j=1}^n$ be such that
\[
\tilde A^{(m)}_{ii}=0\quad\mbox{and}\quad \tilde
A^{(m)}_{ij}=\frac
{A^{(m)}_{ij}-\hat P_{ij}}{\sqrt{(n-1)\hat P_{ij}(1-\hat P_{ij})}},\qquad 1\le i<j\le n. %
\]
\item[(c)] Let $\lambda_1^{(m)}$ and $\lambda_n^{(m)}$ be the largest and
smallest eigenvalues
of $\tilde A^{(m)}$, respectively.
\end{longlist}
\item[4.] Let $(\hat\mu_1,\hat s_1^2)$ and $(\hat\mu_n,\hat s_n^2)$ be the sample
mean and variance of $(\lambda_1^{(m)}:1\le m\le M)$ and $(\lambda
_n^{(m)}:1\le m\le M)$,
respectively.
\item[5.] The bootstrap corrected test statistic is
%
\begin{equation}
\label{eq:boot-stat} T^{(\mathrm{boot})}_{n,K_0}=\mu_{\mathrm{tw}}+s_{\mathrm{tw}}
\max \biggl(\frac{\lambda
_1(\tilde A)-\hat\mu_1}{\hat s_1},-\frac{\lambda_n(\tilde A)-\hat\mu
_n}{\hat s_n} \biggr),
\end{equation}
where $\mu_{\mathrm{tw}}$ and $s_{\mathrm{tw}}$ are the mean and standard
deviation of
the Tracy--Widom distribution.
\end{enumerate}

The fused bootstrap correction is computationally appealing as the
bootstrap sample size
$M$ can be chosen as small as 50. All of our simulations use $M=50$.

\begin{remark*}
The bootstrap correction is based on the empirical observation that although
the finite sample null distribution is different from the theoretical limit,
it has a similar shape, with different location and spread.
Instead of using the theoretical centering and scaling as in \eqref
{eq:A*TW} and
\eqref{eq:TW-A-tilde}, the corresponding bootstrap corrected extreme
eigenvalues are
%
\begin{equation}
\label{eq:boot-eigenvalue} \mu_{\mathrm{tw}}+s_{\mathrm{tw}}\frac{\lambda_1(\tilde A)-\hat\mu_1}{\hat s_1} \quad
\mbox{and}\quad \mu_{\mathrm{tw}}+s_{\mathrm{tw}}\frac{-\lambda_n(\tilde
A)+\hat\mu_n}{\hat s_n}.
\end{equation}
The largest and smallest eigenvalues of $\tilde A$ are individually
corrected using the bootstrap population, because they are individually
shown to have asymptotic Tracy--Widom distribution in Theorem~\ref
{thm:asymp-dist}.
\end{remark*}

%
\begin{figure}[t]

\includegraphics{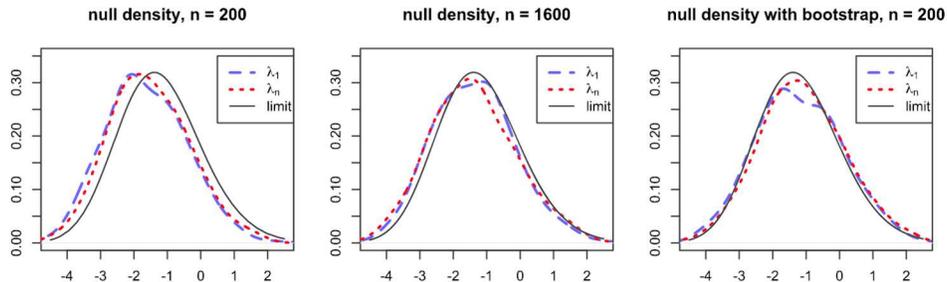}

\caption{The empirical null distributions of scaled and centered
extreme eigenvalues of $\tilde A$ over 1000 repetitions.
Dashed line: largest eigenvalue;
dotted line: smallest eigenvalue;
solid line: theoretical limit distribution.
Left: centered and scaled extreme eigenvalues as in \protect\eqref
{eq:TW-A-tilde} for $n=200$;
middle: centered and scaled extreme eigenvalues as in \protect\eqref
{eq:TW-A-tilde} for $n=1600$;
right: bootstrap corrected extreme eigenvalues as in \protect\eqref
{eq:boot-eigenvalue}.
The stochastic block model used
has two equal-sized communities, and $B_{11}=B_{22}=0.7$, $B_{12}=B_{21}=0.3$.}
\label{fig:null}
\end{figure}

In Figure~\ref{fig:null}, we plot the estimated density of the scaled
and centered extreme eigenvalues of $\tilde A$ calculated from 1000
independent realizations, with and without
bootstrap correction. The stochastic block
model used here has two equal-sized communities, with
$B_{11}=B_{22}=0.7$ and
$B_{12}=B_{21}=0.3$. It is clear that the finite sample null
distribution is systematically
different
from the limiting distribution when $n=200$, and the difference is
reduced but still visible
when $n=1600$. When bootstrap correction is used, the finite sample
null distributions for
both the largest and smallest eigenvalues are
close to the limit even when $n=200$.

\subsection{Simulation 2: Type \textup{I} and type \textup{II} errors}
Now we investigate the type I error of the proposed test under the null
hypothesis
and the power against various alternative distributions.
For each $K_0\in\{2,3,4\}$, we investigate four different models: (i)
the null model, which is a stochastic block model with $K=K_0$ communities;
(ii) a finer stochastic block model (finer SBM) with $K=K_0+1$ communities;
(iii) a degree corrected block model [DCBM, \citet{KarrerN11}] with
$K=K_0$ communities;
and (iv) a mixed membership block model [MMBM, \citet{Airoldi08}] with
$K=K_0$ communities.
For any value of~$K$, the community-wise edge probability matrix $B$ is chosen
such that $B_{kl}=0.2+0.4\times\mathbf1(k=l)$, for all $1\le k,l\le K$.
For the stochastic block model,
the membership vector $g$ is generated by sampling each entry
independently from $\{1,\ldots,K\}$ with equal probability.
For the degree corrected model, the membership vector is generated the same
way as for the stochastic block model, with additional node activeness
parameter $\psi_i$ independently sampled from $\operatorname{Unif}(0,1)$.
In the degree corrected block model, the edge probability between nodes
$i$ and $j$ is $\psi_i\psi_jB_{g_ig_j}$.
For the mixed membership block model, the community mixing probability
$\phi_i$
for each node $i$ is an independent sample from a Dirichlet distribution
with parameter $0.5\times\mathbf e_K$ where $\mathbf e_K$ is a vector
of ones
with length $K$. With such a parameter, each node will tend to favor
one or
two communities so there is a weak community structure. The edge probability
between nodes $i$ and $j$ in the mixed membership block model is $\phi
_i^T B\phi_j$.
For each model, we generate $200$ independent adjacency matrices with $n=1000$
nodes and perform the proposed hypothesis test, with or without
bootstrap correction.
The proportion of rejection at nominal level $0.05$ is summarized in
Table~\ref{tab:level-power}.
%
\begin{table}[t]
\tabcolsep=0pt
\caption{Simulation 2: proportion of rejection at nominal level $0.05$
over $200$ independent samples. The models considered are
\textup{(i)}~Null: the stochastic block model with $K=K_0$ communities;
\textup{(ii)}~Finer SBM: the stochastic block model with $K = K_0+1$ communities;
\textup{(iii)}~DCBM: degree corrected block models with $K=K_0$ communities; and
\textup{(iv)}~MMBM: mixed membership block model with $K=K_0$ communities. The edge probability
between communities $k$ and $l$ is $B_{kl}=0.2+0.4\times\mathbf1 (k=l)$}\label{tab:level-power}
\begin{tabular*}{\tablewidth}{@{\extracolsep{\fill}}@{}lccccc@{}}
\hline
& $\bolds{K_0}$ & \textbf{Null} & \textbf{Finer SBM} & \textbf{DCBM} & \textbf{MMBM}\\
\hline
Without bootstrap & $2$ &$0.02$&$1$&$1$&$1$\phantom{.00}\\
 & $3$ &$0.04$&$1$&$1$&$1$\phantom{.00}\\
& $4$ &$0.03$&$1$&$1$&$0.92$
\\[3pt]
With bootstrap & $2$ &$0.02$&$1$&$1$&$1$\phantom{.00}\\
& $3$ &$0.05$&$1$&$1$&$1$\phantom{.00}\\
& $4$ &$0.06$&$1$&$1$&$0.93$\\
\hline
\end{tabular*}
\end{table}
We observe that the type I error is correctly kept at the nominal
level. The type I error of bootstrap correction method is slightly
closer to the nominal level.
Also we observe that the test can successfully detect all three types
of alternative hypotheses.

\subsection{Simulation 3: Estimating $K$ using sequential testing}
Our third simulation examines the performance of the sequential testing
estimator of $K$
given in \eqref{eq:K-est}.
We use two settings for this simulation.
The first setting concerns different levels of network sparsity, where
the
community-wise connectivity matrices $B$ is given by
$B_{kl}=r (1+2\times\mathbf1(k=l))$.
That is, the edge probability is $3r$ within community and $r$ between
communities.
We consider $r\in\{0.01, 0.02, 0.05, 0.1, 0.2\}$ for different levels
of network
sparsity, and values of $K$ between $2$ and $8$.
For each combination of $K$ and $r$, we generate $200$ independent
adjacency matrices $A$ with $n=1000$ nodes and $K$ equal-sized communities.
The number of communities is estimated for each observation as in \eqref
{eq:K-est} using
threshold $t_n$ corresponding to nominal type I error bound $10^{-4}$.
The proportion
of correct estimates is summarized in Table~\ref{tab:est-1}.
%
\begin{table}[t]
\tabcolsep=0pt
\caption{Simulation 3: proportion of correct estimates of $K$ over 200
simulations under different sparsity levels indexed by $r$. The edge probability
between communities $k$ and $l$ is $r(1+2\times\mathbf1 (k=l))$.
The network size is $n=1000$ with equal sized communities}\label{tab:est-1}
\begin{tabular*}{\tablewidth}{@{\extracolsep{\fill}}@{}ld{1.2}d{1.1}d{1.2}cccd{1.2}d{1.1}cc@{}}
\hline
&\multicolumn{5}{c}{\textbf{With bootstrap}}&\multicolumn{5}{c@{}}{\textbf{Without bootstrap}}\\[-6pt]
&\multicolumn{5}{c}{\hrulefill}&\multicolumn{5}{c@{}}{\hrulefill}\\
$\bolds{r}$& \multicolumn{1}{c}{\textbf{0.01}} & \multicolumn{1}{c}{\textbf{0.02}} & \multicolumn{1}{c}{\textbf{0.05}} &
\multicolumn{1}{c}{\textbf{0.1}} & \multicolumn{1}{c}{\textbf{0.2}} & \multicolumn{1}{c}{\textbf{0.01}} &
\multicolumn{1}{c}{\textbf{0.02}} &\multicolumn{1}{c}{\textbf{0.05}} & \multicolumn{1}{c}{\textbf{0.1}}& \textbf{0.2}\\
\hline
 $K=2$ &1 & 1 & 1 & 1 & 1 & 0.30 & 0.98 & 1 & 1 & 1 \\
 $K=3$ &0.99 & 1 & 1 & 1 & 1 & 0.11 & 0.91 & 1 & 1 & 1 \\
 $K=4$ &0 & 1 & 1 & 1 & 1 & 0.24 & 0.89 & 1 & 1 & 1 \\
 $K=5$ &0 & 0.5 & 1 & 1 & 1 & 0.25 & 0.93 & 1 & 1 & 1 \\
 $K=6$ &0 & 0 & 1 & 1 & 1 & 0.16 & 0.09 & 1 & 1 & 1 \\
 $K=7$ &0 & 0 & 1 & 1 & 1 & 0.04 & 0 & 1 & 1 & 1 \\
 $K=8$ &0 & 0 & 0.71 & 1 & 1 & 0.03 & 0 & 0.9 & 1 & 1 \\
\hline
\end{tabular*}
\end{table}
The sequential testing estimator with bootstrap correction works well
for $K=2,3$ at all sparsity levels.
When $K$ gets larger, both methods require denser models
to consistently estimate $K$. When the model is moderately dense, both
methods work well for all values of $K$.
For very sparse models, the null distribution without bootstrap is
biased and the sequential testing method tends to pick larger values of~$K$.

In the second setting, the focus is on different types of block structures.
To this end,
for each $K\in\{2,3,4\}$ we generate matrices $B$ whose diagonal and
upper diagonal entries are
independently drawn from a uniform distribution between $0$ and $0.5$.
The success of spectral clustering requires the smallest singular value
of $B$ to be bounded away
from zero, so we only use those $B$ matrices whose smallest singular
values are at least $0.1$.
The membership vector $g$ is generated by sampling each entry
independently from $\{1,\ldots,K\}$ with equal probability.
For each $K$ and network size $n=500$ and $n=1000$,
we generate $200$ independent adjacency matrices using random $B$ and
$g$ described above.
Similarly, $K$ is estimated as in \eqref{eq:K-est} using
threshold $t_n$ corresponding to nominal type I error bound $10^{-4}$.
In Table~\ref{tab:est-K}, we summarize the proportion of correct estimates.
The proposed test can correctly estimate the number of communities in a
very large
proportion of these randomly generated models.
In general,
the bootstrap correction helps improve the estimation accuracy.

%
\begin{table}[b]
\tabcolsep=10pt
\caption{Simulation 3: proportion of correct estimates of $K$ over 200 simulations
with randomly generated matrices $B$ and membership vectors $g$}\label{tab:est-K}
\begin{tabular*}{\tablewidth}{@{\extracolsep{\fill}}@{}lcccccc@{}}
\hline
&\multicolumn{3}{c}{\textbf{With bootstrap}} & \multicolumn{3}{c@{}}{\textbf{Without bootstrap}}\\[-6pt]
&\multicolumn{3}{c}{\hrulefill} & \multicolumn{3}{c@{}}{\hrulefill}\\
$\bolds{K}$& \textbf{2} & \textbf{3} & \textbf{4} & \textbf{2} & \textbf{3} & \textbf{4}\\
\hline
$n=500$ &$0.99$&$0.90$&$0.76$&$0.91$&$0.84$&$0.74$\\
$n=1000$ &$1$\phantom{.00}&$1$\phantom{.00}&$0.93$&$0.98$&$0.93$&$0.90$\\
\hline
\end{tabular*}
\end{table}

\subsection{The political blog data}
The political blog data \citep{PolBlog} records hyperlinks between web
blogs shortly before the 2004 US presidential election. It has been
used widely in the network community detection literature
as an example of significant within-community node degree variation
[see \citet{KarrerN11,ZhaoLZ12,Jin12}, e.g.]. It is widely believed that
a degree corrected block model is more suitable for this data, rather
than a regular stochastic block model. \citet{Yan14} used a likelihood
ratio method to choose the degree corrected model over the regular
stochastic block model. Theoretical justification of the $\chi^2$
approximation used in this method is still an open problem, and
maximizing the likelihood is computationally demanding. Following
common practice, we consider the largest connected component of the
political blog data. There are 1222 nodes with community sizes 586 and
636. We set $\hat g$ to be the true labeling given in the data---the
results are similar for $\hat g$ estimated from the data.
Under the null hypothesis that the data is generated from a stochastic
block model of two communities, the test statistic is $1172.3$ for the
original test and $491.5$ for the bootstrap corrected test,
both indicating strong evidence to reject the null hypothesis. In
addition, we apply the sequential testing procedure at type I error
level $10^{-5}$, with block model parameters estimated by spectral
clustering using two leading eigenvectors of the adjacency matrix. The
procedure partitions the nodes into 17 groups. $\mathrm{Sixteen}$ of
these estimated groups mostly contain nodes from one true community,
with $8$ groups for each community and stratified by degrees.
The additional estimated group contains nodes with very small degrees,
whose community memberships are very hard to recover.

\section{Discussion}
\label{sec:discussion}
The goodness-of-fit test developed in this paper is an attempt to
perform principled
statistical inference for stochastic block models. The test statistic
reflects a fundamental
difference between network models and traditional statistical models on
independent
individuals. In traditional independent and identically distributed
data samples, the
goodness-of-fit is usually assessed by the sum of residuals or squared
residuals.
For stochastic block models, the residual is a matrix, where the signal
is not carried in
the sum of individual residuals but is determined by how these
residuals align across the rows
and columns. For example, suppose $A$ is generated from a stochastic
block model with two
communities and we want to test if $K=1$. If we simply treat the upper
diagonal entries of $A$ as
independent Bernoulli variables, the goodness-of-fit test reduces to
testing whether
the $n(n-1)/2$ upper diagonal entries look like an independent sample
of a Bernoulli variable.
Such tests have
little power in detecting the block structure. On the other hand, the
extreme singular value
of the residual matrix accurately captures the block structure.
This is an example of detecting low-rank mean effect from a noisy
random matrix
using its extreme eigenvalues. Other examples using the similar idea include
\citet{Kargin14} for reduced
rank multivariate regression and
\citet{Montanari14} for the Gaussian hidden clique problem.
It would be interesting to further develop goodness-of-fit testing
methods for more realistic null hypotheses, such as the degree
corrected block model or even the nonparametric graphon model \citep{WolfeO13}.

It is possible to extend the method and theory developed in this paper to
certain sparse stochastic block models.
Consider sparse
stochastic block models with $B=\rho_nB_0$ where
the entries of $B_0$ are of order $1$ and $\rho_n\downarrow0$
controls the overall network sparsity.
Most random matrix theory used in this paper
(namely, Lemmas~\ref{lem:TW-A*}, \ref{lem:eigenvalue-count}, \ref{lem:deviation}) has been
developed for moderately sparse stochastic block models with $\rho_n
\gg n^{-1/3}$
in
\citeauthor{ER-RMT-1} (\citeyear{ER-RMT-1,ER-RMT-2}).
However, existing arguments do not guarantee isotropic delocalization
of eigenvectors (Lemma~\ref{lem:deloc})
due to the heavy tail of the normalized adjacency matrix entries
$(A_{ij}-P_{ij})/[P_{ij}(1-P_{ij})]$.
The possibility of proving such a result using modified techniques has
been mentioned in
\citet{General}.


\begin{appendix}\label{sec:proof}
\section*{Appendix: Proofs}
\subsection*{Additional notation}
Let $(\lambda_j^*,u_j^*)_{j=1}^n$ be the eigenvalue-eigenvector pairs
of~$\tilde A^*$ such that $\lambda_1^*\ge\lambda_2^*\ge\cdots\ge\lambda_n^*$.
For a pair of random sequences $(a_n)$ and~$(b_n)$, we write
$a_n=\tilde O_P(b_n)$
if for any $\varepsilon>0$ and $D>0$ there exists $n_0=n_0(\varepsilon,D)$
such that
\[
P\bigl(a_n\ge n^\varepsilon b_n \bigr)\le
n^{-D}\qquad\mbox{for all } n\ge n_0. %
\]
For any matrix $M$ with singular value decomposition $M=\sum_{j} \sigma
_j u_j v_j^T$, define $\llvert  M\rrvert  =\sum_{j}\llvert  \sigma_j\rrvert  u_j v_j^T$. We will use
$c$ and $C$ to denote
positive constants independent of $n$, which may vary from line to line.

\subsection{Results from random matrix theory}\label{sec:RMT}
We first collect some useful results from random matrix theory
regarding the
distributions of the eigenvalues and eigenvectors of $\tilde A^*$.

\begin{lemma}[{[Asymptotic distributions of $\lambda_1(\tilde A^*)$ and $\lambda_n(\tilde A^*)$]}]\label{lem:TW-A*} For $\tilde A^*$ defined in
\eqref{eq:A*} we have
\[
n^{2/3}\bigl(\lambda_1\bigl(\tilde A^*\bigr)-2\bigr)
\rightsquigarrow TW_1,\qquad n^{2/3}\bigl(-
\lambda_n\bigl(\tilde A^*\bigr)-2\bigr)\rightsquigarrow
TW_1.
\]
\end{lemma}

\begin{pf}
Let $G^*$ be an $n\times n$ symmetric matrix
whose upper diagonal entries are independent
normal with mean zero and variance $1/(n-1)$, and
diagonal entries are zero. Then $\tilde A^*$ and
$G^*$ have the same first and second moments. According to
Theorem 2.4 of \citet{Rigidity}, we know that
$n^{2/3}(\lambda_1(\tilde A^*)-2)$ and $n^{2/3}(\lambda_1(G^*)-2)$
have the same limiting distribution.
But $n^{2/3}(\lambda_1(G^*)-2)\rightsquigarrow TW_1$ according to
\citet{LeeY14}. The same argument applies to $\lambda_n(\tilde A^*)$.
\end{pf}

%
\begin{lemma}[(Eigenvector delocalization)]\label{lem:deloc}
For each deterministic unit vector $u$ and each $1\le j \le n$, for
any $\varepsilon>0$ and $D>0$ there exists $n_0=n_0(\varepsilon,D)$ such that
\[
P \bigl[\bigl(u^T u_j^*\bigr)^2\ge
n^{-1+\varepsilon} \bigr] \le n^{-D}\qquad\mbox {for all } n\ge
n_0. %
\]
\end{lemma}

It is worth noting that the above result is uniform over $j$ and $u$ in
the sense that $n_0(\varepsilon,D)$
does\vspace*{1pt} not depend on $u$ or $j$. Lemma~\ref{lem:deloc} can be
equivalently stated as
$(u^T u_j^*)^2=\tilde O_P(n^{-1})$ uniformly over all $u_j^*$ ($1\le
j\le n$) and
all deterministic unit vector $u$.

Lemma~\ref{lem:deloc} is Theorem 2.16 of \citet{Isotropic}.
Although \citet{Isotropic} requires the diagonal entries of $\tilde A^*$
to have positive variance, their Theorem 2.16 is a consequence of the
local semicircle law [Theorem 2.12 of \citet{Isotropic}], which can be
established
for matrices with zero diagonals using the result of \citet{General}.
See also the discussion
in \citet{BickelS13test}.

%
\begin{lemma}[(Counting large eigenvalues)]\label{lem:eigenvalue-count}
Let $c_n$ be a possibly random number of order $o_P(1)$ and $m(c_n)$
be the number of eigenvalues of $\tilde
A^*$
larger than $\lambda_1^*-c_n$.
Then $m(c_n)=O_P(n c_n^{3/2}) +\tilde O_P(1)$.
\end{lemma}

Lemma~\ref{lem:eigenvalue-count} extends equation (26) of \citet{BickelS13test}.

\begin{pf*}{Proof of Lemma \ref{lem:eigenvalue-count}}
For any $a<b<5$, let $N^*(a,b)$ be the number of
eigenvalues of $\tilde A^*$ in the interval $(a,b]$,
and $N(a,b)=n\int_{a}^b \rho_{\mathrm{sc}}(x)\,dx$
where $\rho_{\mathrm{sc}}(x)=(1/2\pi)((4-x^2)_+)^{1/2}$ is the density of
the semicircle law. Let $\delta(a,b)=N^*(a,b)-N(a,b)$ then according to
Theorem 2.2 of \citet{Rigidity} we have $\sup_{a,b<5} \llvert  \delta
(a,b)\rrvert  =\tilde O_P(1)$.
Then, conditioning on the event that $\{\llvert  2-\lambda_1^*\rrvert  +c_n\le1\}$,
we have
\begin{eqnarray*}
m(c_n) &=& N^*\bigl(\lambda_1^*-c_n,
\lambda_1^*\bigr)
\\
&=& N\bigl(\lambda_1^*-c_n,
\lambda_1^*\bigr) + \sup_{a,b<5}\bigl\llvert
\delta(a,b)\bigr\rrvert
\\
&\le& n\int_{2-(2-\lambda_1^*)-c_n}^2 \bigl(\bigl(4-x^2
\bigr)_+\bigr)^{1/2}\,dx + \tilde O_P(1)
\\
&\le& 2n \bigl(c_n+\bigl\llvert 2-\lambda_1^*\bigr
\rrvert \bigr)^{3/2}+\tilde O_P(1)
\\
&\le& O\bigl(n
c_n^{3/2}\bigr) +\tilde O_P(1). 
\end{eqnarray*}\upqed
\end{pf*}

The claimed result follows by observing that the event $\{\llvert  2-\lambda
_1^*\rrvert  +c_n\le1\}$
has probability $1-o(1)$.

%
\begin{lemma}[(Deviation of largest singular value)]\label{lem:deviation}
There exists absolute positive constants $a$, $b$, $c$, $C$, such that
\[
P \bigl[n^{2/3}\bigl(\sigma_1\bigl(\tilde A^*\bigr)-2
\bigr)\ge(\log n)^{a\log\log n} \bigr] \le C \exp \bigl[-b(\log n)^{c\log\log n}
\bigr].
\]
\end{lemma}

Lemma~\ref{lem:deviation} is a direct consequence of equation (2.22) in
\citet{Rigidity}.
We can simplify the statement so that there exists an
absolute constant $b>0$ such that for any $\varepsilon>0$
%
\begin{equation}
\label{eq:deviation} P \bigl[n^{2/3}\bigl(\sigma_1\bigl(\tilde
A^*\bigr)-2\bigr)\ge n^\varepsilon \bigr] = O\bigl(n^{-b}\bigr).
\end{equation}

\subsection{Proof of asymptotic null distribution}\label{sec:proof-null}
Now we provide proofs for theoretical results in Section~\ref
{sec:theory}. Here,
we omit the dependence on $n$ in $g$, $B$ and $K$ for simplicity.

\begin{pf*}{Proof of Theorem~\ref{thm:asymp-dist}}
The consistency of $\hat g$ allows us to focus on the event $\hat g=g$.

We will prove the claim for $\lambda_1(\tilde A)$. The other claim can
be proved by applying the same argument on $-\tilde A$.

Let $\tilde A'\in\mathbb R^{n\times n}$ be such that
%
\begin{equation}
\label{eq:tilde-A'} \tilde A'_{ij} = \cases{ \displaystyle
\frac{A_{ij}-\hat P_{ij}}{\sqrt{(n-1)P_{ij}(1-P_{ij})}},&\quad$i\neq j$,
\cr
\displaystyle\frac{P_{ii}-\hat P_{ii}}{\sqrt{(n-1)P_{ii}(1-P_{ii})}},&
\quad$i=j$.}
\end{equation}
Thus, $\tilde A'=\tilde A^*+\Delta'$, where
$\Delta'_{ij}=(P_{ij}-\hat P_{ij})/\sqrt{(n-1)P_{ij}(1-P_{ij})}$.
Because $\Delta'$ is a $K\times K$ block-wise constant symmetric
matrix, its rank is at most $K$, and the corresponding principal subspace
is spanned by $(\theta_1,\ldots,\theta_K)$, where
$\theta_k\in\mathbb R^{n}$ is the unit norm indicator of the $k$th
community in $g$.
That is, the $i$th entry of $\theta_k$ is $n_k^{-1/2}$ if $g_i=k$ and
zero otherwise,
where $n_k$ is the size of the $k$th community.

The consistency of $\hat g$ implies that with probability tending to
one, for each $1\le k,k'\le K$,
$\hat B_{k,k'}$ is the sample mean of independent Bernoulli random
variables with parameter $B_{k,k'}$ and
sample size of order $(n/K)^2$. Thus, standard large deviation
inequalities such as Bernstein's inequality or Hoeffding's inequality
suggest that $\sup_{k,k'}\llvert  \hat B_{k,k'}-B_{k,k'}\rrvert  =o_P(K \log n / n)$,
which implies that
$\sup_{i,j}\llvert  \hat P_{ij}-P_{ij}\rrvert  =o_P(K\log n/n)$.
Note here the $o_P$ statement goes through a union bound over $K^2$
terms, which is
valid since the tail probability bound for $\hat P_{ij}-P_{ij}$ can be
made exponentially small in $n$.
Let $\Delta' = \Theta\Gamma\Theta^T$, where
$\Theta=(\theta_1,\ldots,\theta_K)$ and $\Gamma$ is a $K\times K$
symmetric matrix.
Then each entry of $\Gamma$ is $o_P(n^{-1/2}\log n)$, and hence $\llVert
\Gamma\rrVert  =o_P(Kn^{-1/2}\log n)$.

We will show that
%
\begin{equation}
\label{eq:A'-A*} \lambda_1\bigl(\tilde A'\bigr)=
\lambda_1\bigl(\tilde A^*\bigr)+o_P\bigl(n^{-2/3}
\bigr),
\end{equation}
by establishing a lower and upper bound on $\lambda_1(\tilde A')$. Both
parts uses
the eigenvector delocalization result (Lemma~\ref{lem:deloc}) as follows.
Let $\Theta= (\theta_1,\ldots,\theta_K)$, then, uniformly over $j$ we have
%
\begin{equation}
\label{eq:theta^Tu^*} \bigl\llVert \Theta^T u_j^*\bigr\rrVert
_2^2=\sum_{k=1}^K
\bigl(\theta_k^T u_j^*\bigr)^2 =
\tilde O_P\bigl(K n^{-1}\bigr),
\end{equation}
and hence
%
\begin{eqnarray}\label{eq:delta-u^*}
\bigl\llvert \bigl(u_j^*\bigr)^T
\Delta' u_j^*\bigr\rrvert &\le& \bigl\llvert \bigl(
\Theta^T u_j^*\bigr)^T \Gamma\bigl(
\Theta^T u_j^*\bigr)\bigr\rrvert \le\bigl\llVert
\Theta^T u_j^*\bigr\rrVert _2^2
\llVert \Gamma\rrVert
\nonumber\\[-8pt]\\[-8pt]\nonumber
&=& \tilde O_P\bigl(K^2 n^{-3/2}\log n\bigr).
\end{eqnarray}
Here, the $\tilde O_P$ statement in \eqref{eq:theta^Tu^*}
holds when taking union bound over $K$ terms by choosing $D$ large
enough in Lemma~\ref{lem:deloc}.

First, we provide a lower bound on $\lambda_1(\tilde A')$:
%
\begin{eqnarray}\label{eq:A'-A*-l}
\lambda_1\bigl(\tilde A'\bigr)&\ge&
\bigl(u_1^*\bigr)^T \tilde A'
u_1^* = \lambda_1^* + \bigl(u_1^*
\bigr)^T \Delta' u_1^*
\nonumber
\\
&\ge& \lambda_1^*-\tilde
O_P\bigl(K^2 n^{-3/2}\log n\bigr)
\\
&\ge&\lambda_1^*-o_P\bigl(n^{-2/3}\bigr),\nonumber
\end{eqnarray}
where the last inequality uses the assumed upper bound on the rate at
which $K$ grows with $n$, and
the second last inequality uses \eqref{eq:delta-u^*}.

Next, we provide an upper bound of $\lambda_1(\tilde A')$.
For any unit vector $u\in\mathbb R^{n}$, let $(a_1,\ldots,a_n)$
be a unit vector in $\mathbb R^{n}$ such that
\[
u=\sum_{j=1}^n a_j
u_j^*.
\]

Let $m$ be the number of $\lambda_j^*$'s
in the interval $(\lambda_1^*-2\llVert  \Delta'\rrVert,\lambda_1^*]$,
and
$u_1=\sum_{j=1}^m a_j u_j^*$, $u_2=\sum_{j=m+1}^n a_j u_j^*$.
Then
%
\begin{eqnarray}\label{eq:A'-A*-u}
u^T \tilde A' u &=& u^T\tilde A^* u +
u^T \Delta' u
\nonumber
\\
&\le& \lambda_1^*\sum_{j=1}^m
a_j^2 + \bigl(\lambda_{1}^*-2\bigl\llVert
\Delta'\bigr\rrVert \bigr)\sum_{j=m+1}^n
a_j^2 + 2 u_1^T \bigl\llvert
\Delta'\bigr\rrvert u_1 + 2 u_2^T
\bigl\llvert \Delta'\bigr\rrvert u_2
\nonumber
\\
&\le& \lambda_1^*\sum_{j=1}^m
a_j^2 + \bigl(\lambda_{1}^*-2\bigl\llVert
\Delta'\bigr\rrVert \bigr)\sum_{j=m+1}^n
a_j^2
\nonumber
\\
&&{} + 2m \sum_{j=1}^m
a_j^2\bigl(u_j^*\bigr)^T \bigl
\llvert \Delta'\bigr\rrvert u_j^* + 2
u_2^T \bigl\llvert \Delta'\bigr\rrvert
u_2
\nonumber
\\
&\le&\lambda_1^*\sum_{j=1}^m
a_j^2 + \bigl(\lambda_{1}^*-2\bigl\llVert
\Delta'\bigr\rrVert \bigr)\sum_{j=m+1}^n
a_j^2
\\
&&{}+ 2 m \tilde O_P\bigl(K^2 n^{-3/2}\log n
\bigr)\sum_{j=1}^m a_j^2
+ 2 \bigl\llVert \Delta'\bigr\rrVert \sum
_{j=m+1}^n a_j^2
\nonumber
\\
&\le&\lambda_1^*+m\tilde O_P\bigl(K^2
n^{-3/2}\log n\bigr)
\nonumber
\\
&\le& \lambda_1^*+ \bigl(O\bigl(n\bigl\llVert \Delta'
\bigr\rrVert ^{3/2}\bigr)+\tilde O_P(1) \bigr)\tilde
O_P\bigl(K^2 n^{-3/2}\log n\bigr)
\nonumber
\\
&=& \lambda_1^*+\tilde O_P\bigl(K^{7/2}(\log
n)^{5/2} n^{-5/4}\bigr),\nonumber
\end{eqnarray}
where the\vspace*{2pt} third inequality uses \eqref{eq:delta-u^*} and uniformity
over $j$, and the
second last line uses Lemma~\ref{lem:eigenvalue-count} together with $\llVert
\Delta'\rrVert  =o_P(K n^{-1/2}\log n)$.

Thus, \eqref{eq:A'-A*} is established by combining \eqref{eq:A'-A*-l}
and \eqref{eq:A'-A*-u}, provided that $K=O(n^{1/6-\tau})$ for some\vspace*{1pt}
small positive $\tau$.\vspace*{1pt}

Next,\vspace*{1pt} we show that $\lambda_1(\tilde A)=\lambda_1(\tilde A')+o_P(n^{-2/3})$.
Let $\tilde A''=\tilde A'-\operatorname{diag}(\tilde A')$.
Consider the block representation of $\tilde A$:
\[
\tilde A=(\tilde A_{(k,l)})_{k,l=1}^K,
\]
where $\tilde A_{(k,l)}$ is the submatrix corresponding to the rows
in community $k$ and columns in community $l$. Similar block representations
can be defined for $\tilde A''$.
It is obvious that
\[
\tilde A_{(k,l)}=\tilde A''_{(k,l)}
\frac{\sqrt{B_{kl}(1-B_{kl})}}{\sqrt
{\hat B_{kl}(1-\hat B_{kl})}}=\tilde A''_{(k,l)}
\bigl(1+o_P\bigl(Kn^{-1}\log n\bigr)\bigr). %
\]
Therefore,
\begin{eqnarray*}
\bigl\llVert \tilde A-\tilde A''\bigr\rrVert &\le& K
\max_{k,l} \bigl\llVert \tilde A_{(k,l)}-\tilde
A''_{(k,l)}\bigr\rrVert \le o_{P}
\bigl(Kn^{-1}\log n\bigr) K \sum_{k,l}\bigl
\llVert \tilde A''_{(k,l)}\bigr\rrVert
\\
&\le&o_P\bigl(K^2 n^{-1}\log n\bigr) \bigl
\llVert \tilde A''\bigr\rrVert \le o_P
\bigl(K^2 n^{-1}\log n\bigr) \bigl(\bigl\llVert \tilde
A'\bigr\rrVert +\bigl\llVert \operatorname{diag}\bigl(\tilde A'
\bigr)\bigr\rrVert \bigr)
\\
&\le&o_P\bigl(K^2 n^{-1}\log n\bigr)
\bigl(O_P(1)+O_P\bigl(K n^{-3/2}\log n\bigr)
\bigr)
\\
&=& o_P\bigl(K^2 n^{-1}\log n
\bigr)=o_P\bigl(n^{-2/3}\bigr).
\end{eqnarray*}
%

Then
%
\begin{equation}
\bigl\llVert \tilde A-\tilde A'\bigr\rrVert \le\bigl\llVert \tilde
A-\tilde A''\bigr\rrVert +\bigl\llVert \operatorname{diag}\bigl(
\tilde A'\bigr)\bigr\rrVert =o_P\bigl(n^{-2/3}
\bigr).\label{eq:A-A'}
\end{equation}
Combining \eqref{eq:A'-A*} and \eqref{eq:A-A'}, we have
%
\begin{equation}
\label{eq:A-A*} \lambda_1(\tilde A) = \lambda_1\bigl(
\tilde A^*\bigr) + o_P\bigl(n^{-2/3}\bigr).
\end{equation}
Now applying Lemma~\ref{lem:TW-A*} and combining with \eqref{eq:A-A*}
we have
\begin{eqnarray*}
n^{2/3}\bigl(\lambda_1(\tilde A)-2\bigr) &
\rightsquigarrow & TW_1.
\end{eqnarray*}
\upqed
\end{pf*}

\subsection{Proof of power and consistency}\label{sec:proof-power}

\mbox{}

\begin{pf*}{Proof of Theorem~\ref{thm:consistency}}
For all $1\le l\le K$, $1\le k\le K_0$, let $\mathcal N_l=\{i: g_i=l\}$,
$\hat{\mathcal N}_k=\{i:\hat g_i=k\}$
and $\hat{\mathcal N}_{k,l}=\{i:\hat g_i=k,g_i=l\}$.
For each $1\le l\le K$, $\mathcal N_l$ is partitioned
into $\{\hat{\mathcal N}_{k,l}:1\le k\le K_0\}$. Thus,
for each $1\le l\le K$ there exists a $k_l$ such that $1\le k_l \le
K_0$ and $\llvert  \hat{\mathcal N}_{k_l,l}\rrvert  \ge\llvert  \mathcal N_l\rrvert  /K_0\ge c_0 n /
(K\times K_0)\ge c_0 n K^{-2}$.
Because $K_0<K$, there exist $l_1$ and $l_2$ such that
$k_{l_1}=k_{l_2}=k$.
Since $B\in\mathcal B_K$, there exists an $l_3$ such that
$B_{l_1, l_3}\neq B_{l_2,l_3}$. Let $k'=k_{l_3}$ and we have $\llvert  \hat
{\mathcal N}_{k',l_3}\rrvert  \ge c_0 n K^{-2}$.

Let $\tilde A^{(0)}$ be the submatrix of $\tilde A$ consisting the rows in
$\hat{\mathcal N}_{k,l_1}\cup\hat{\mathcal N}_{k,l_2}$,
and the columns in $\hat{\mathcal N}_{k',l_3}$. Define $A^{(0)}$, $\hat
P^{(0)}$, and $P^{(0)}$
correspondingly.

When\vspace*{1pt} $k\neq k'$, or $k=k'$ but $l_3\notin\{l_1,l_2\}$, the submatrix
$A^{(0)}$ contains only off-diagonal entries of $A$.
Therefore, $\hat P^{(0)}$ is a constant matrix in that all of its
entries are equal.
We have
%
\begin{eqnarray}\label{eq:alt-lower-bound}
\llVert \tilde A\rrVert &\ge& \bigl\llVert \tilde A^{(0)}\bigr\rrVert
\ge n^{-1/2}\bigl\llVert A^{(0)}-\hat P^{(0)}\bigr
\rrVert
\nonumber
\\
&\ge& n^{-1/2} \bigl(\bigl\llVert P^{(0)}-\hat P^{(0)}
\bigr\rrVert -\bigl\llVert A^{(0)}-P^{(0)}\bigr\rrVert \bigr)
\nonumber\\[-8pt]\\[-8pt]\nonumber
&\ge& n^{-1/2} \bigl(\bigl\llVert P^{(0)}-\hat P^{(0)}
\bigr\rrVert -O_P\bigl(n^{1/2}\bigr) \bigr)
\\
&\ge& n^{-1/2}\bigl(\delta_B c_0 nK^{-2}/2
-O_P\bigl(n^{1/2}\bigr)\bigr).\nonumber
\end{eqnarray}
To obtain the last inequality, first note that $P^{(0)}$ has
two distinct blocks each with at least $c_0 nK^{-2}$ rows and at least
$c_0 nK^{-2}$ columns. Each of these two blocks has constant entries
and at least one of them
has absolute entry value at least $\delta_B/2$. Thus,
$\llVert  P^{(0)}-\hat P^{(0)}\rrVert  \ge\delta_B c_0 n K^{-2} / 2$.

When $k=k'$ and $l_3\in\{l_1,l_2\}$,
the submatrix $A^{(0)}$ defined above contains diagonal entries of $A$.
The\vspace*{1pt} corresponding entries of $\hat P^{(0)}$ are zero. These
zero entries causes an additional $O(1)$ term in $\llVert  \hat
P^{(0)}-P^{(0)}\rrVert  $
and \eqref{eq:alt-lower-bound} still goes through.
\end{pf*}

\begin{pf*}{Proof of Corollary~\ref{cor:consistency}}
Following the notation in the proof of Theorem~\ref{thm:consistency},
for any $K_0 < K$
we have, in view of \eqref{eq:alt-lower-bound} and letting
$C=\inf_n \delta_B c_0/2$,
%
\begin{eqnarray}\label{eq:consist-1}
P(T_{n,K_0} < t_n) &= & P \bigl[n^{2/3}\bigl(\llVert
\tilde A\rrVert -2\bigr)<t_n \bigr]= P \bigl[\llVert \tilde A\rrVert
<n^{-2/3}t_n+2 \bigr]
\nonumber
\\
&\le& P \bigl[n^{-1/2} \bigl( \bigl\llVert P^{(0)}-\hat
P^{(0)}\bigr\rrVert -\bigl\llVert A^{(0)}-P^{(0)}\bigr
\rrVert \bigr)\le n^{-2/3}t_n+2 \bigr]\nonumber
\\
&\le& P \bigl[n^{-1/2} \bigl\llVert A^{(0)}-P^{(0)}
\bigr\rrVert \ge Cn^{1/2}K^{-2}-n^{-2/3}t_n-2
\bigr]\nonumber
\\
&\le& n^{-1},\nonumber
\end{eqnarray}
the\vspace*{1pt} last inequality is obtained by first using the assumption
$K=O(n^{1/6-\tau})$, and $t_n=O(n^{5/6})$,
so that $Cn^{1/2}K^{-2}+n^{-2/3}t_n+2\ge n^{1/6}$ for large $n$, and
then applying operator norm deviation bound results such as Theorem 5.2
of \citet{LeiR14} [see also Theorem 3.4 of \citet{Chatterjee15}].

Therefore,
\[
P(\hat K < K)\le\sum_{K_0=1}^{K-1}P(T_{n,K_0}<t_n)
\le n^{-1}(K-1)=o(1).
\]

On the other hand,
\begin{eqnarray*}
P(\hat K>K)&\le& P(T_{n,K}\ge t_n) = P\bigl(n^{2/3}
\bigl(\sigma_1(\tilde A)-2\bigr)\ge t_n\bigr)
\\
&\le& P\bigl(n^{2/3}\bigl(\sigma_1\bigl(\tilde A^*\bigr)-2
\bigr)\ge t_n/2\bigr)+ P\bigl(n^{2/3}\bigl\llvert
\sigma_1\bigl(\tilde A^*\bigr)-\sigma_1(\tilde A)\bigr
\rrvert \ge t_n/2\bigr)
\\
&=&o(1),
\end{eqnarray*}
where the first probability is controlled using Lemma~\ref{lem:deviation}
and the second probability is controlled using \eqref{eq:A-A*}
and its analogous result for $\lambda_n(\tilde A)-\lambda_n(\tilde A^*)$.
\end{pf*}

\subsection{Asymptotic power against degree corrected block models}\label{sec:proof-dcbm}

\mbox{}

\begin{pf*}{Proof of Theorem~\ref{thm:dcbm}}
Recall that $\hat{\mathcal N}_{l,k}=\{i:g_i=k,\hat g_i=l\}$ ($1\le
l\le K_0$, $1\le k\le K$).
Let $\tilde\phi_{k,\hat{\mathcal N}_{l,k}}$ be the subvector of $\tilde
\phi_{k}$ on the entries
in $\hat{\mathcal N}_{l,k}$.

Let
$\eta_l=\tilde\phi_{k^*,\hat{\mathcal N}_{l,k^*}}$ for each $1\le l\le K_0$.
By definition of $\mathcal E$, $\sum_{l=1}^{K_0}\mathcal E(\eta
_l,1)\ge\mathcal E(\tilde\psi_{k^*},K_0)$, and hence there
exists an $l^*$ such that $\mathcal E(\eta_{l^*},1)\ge\mathcal
E(\tilde\phi_{k^*},K_0)/K_0$.\vspace*{1pt}

For simplicity, denote $\eta=\eta_{l^*}$ and $\bar\eta= \eta/\llVert  \eta\rrVert  $.
Let $m=\llvert  \hat{\mathcal N}_{l^*,k^*}\rrvert  $ and define $\mathbf e_m$ as the
$1\times m$ vector with $1/\sqrt{m}$
on each entry.
Therefore, we have
%
\begin{equation}
\label{eq:dcbm-pf-1} \llVert \bar\eta-\mathbf e_m\rrVert ^2 \ge
\mathcal E(\bar\eta,1) = \llVert \eta\rrVert ^{-2}\mathcal E(\eta,1) \ge
\llVert \eta\rrVert ^{-2}\mathcal E(\tilde\phi_{k^*},K_0)/K_0.
\end{equation}
Because $\bar\eta$ and $\mathbf e_m$ both have unit $\ell_2$ norm,
\eqref{eq:dcbm-pf-1}
implies that
\[
\bigl\llvert \mathbf e_m^T\bar\eta\bigr\rrvert \le1-
\frac{\mathcal E(\tilde\phi
_{k^*},K_0)}{2K_0\llVert  \eta\rrVert  ^2}. %
\]
Let $u = (\bar\eta- \mathbf e_m \mathbf e_m^T \bar\eta)/\llVert  \bar\eta-
\mathbf e_m \mathbf e_m^T \bar\eta\rrVert  $, then
%
\begin{equation}
\label{eq:dcbm-pf-2} u^T \eta= u^T \bar\eta\llVert \eta\rrVert
=\llVert \eta\rrVert \frac{1-(\mathbf e_m^T \bar
\eta)^2}{\llVert  \bar\eta- \mathbf e_m \mathbf e_m^T \bar\eta\rrVert  }\ge\frac
{\mathcal E(\tilde\phi_{k^*},K_0)}{2K_0\llVert  \eta\rrVert  } \ge
\frac{\mathcal E(\tilde\phi_{k^*},K_0)}{2K_0}.
\end{equation}

Now let $k'$ be such that $B_{k^*,k'}=\llVert  B_{k^*,\cdot}\rrVert  _{\infty}$.
There exists an $l'$ such that $\llVert  \tilde\phi_{k',\hat{\mathcal
N}_{l',k'}}\rrVert  \ge K_0^{-1/2}$.

Let $A^{(0)}$ be the submatrix of $A$ corresponding to the rows in
$\hat{\mathcal N}_{l^*,k^*}$ and columns in $\hat{\mathcal N}_{l',k'}$,
and define $\tilde A^{(0)}$, $P^{(0)}$, $\hat P^{(0)}$ similarly.
Thus, by construction we have, letting $m'=\llvert  \hat{\mathcal N}_{l',k'}\rrvert  $,
\[
P^{(0)}= \llVert \phi_{k^*}\rrVert \llVert
\phi_{k'}\rrVert B_{k^*k'} \eta\tilde\phi_{k',\hat
{\mathcal N}_{l',k'}}^T,\qquad \hat P^{(0)} = \hat B_{l^*l'}\sqrt{m m'}
\mathbf e_m \mathbf e_{m'}^T. %
\]

Observing that $u^T\mathbf e_m=0$, we have
\begin{eqnarray*}
\bigl\llVert u^T\bigl(P^{(0)}-\hat P^{(0)}\bigr)
\bigr\rrVert &= & \bigl\llVert u^T P^{(0)}\bigr\rrVert =
\llVert \psi_{k^*}\rrVert \llVert \psi_{k'}\rrVert
B_{k^*k'} \bigl\llvert u^T \eta\bigr\rrvert \llVert \tilde
\phi_{k',\hat{\mathcal
N}_{k',l'}}\rrVert
\\
&\ge& \llVert \psi_{k^*}\rrVert \llVert \psi_{k'}\rrVert
B_{k^*k'} \frac{\mathcal E(\tilde
\phi_{k^*},K_0)}{2K_0^{3/2}}
\\
&\ge& \kappa_n n \llVert B_{k^*,\cdot}\rrVert _{\infty}
\frac{\mathcal E(\tilde\phi
_{k^*},K_0)}{2K_0^{3/2}}.
\end{eqnarray*}
The claimed result follows by observing that
\begin{eqnarray*}
\bigl\llVert \tilde A^{(0)}\bigr\rrVert &\ge& n^{-1/2} \bigl(
\bigl\llVert P^{(0)}-\hat P^{(0)}\bigr\rrVert -\bigl\llVert
A^{(0)}-P^{(0)}\bigr\rrVert \bigr)
\end{eqnarray*}
and $\llVert  A^{(0)}-P^{(0)}\rrVert  =O_P(\sqrt{n})$.
\end{pf*}
\end{appendix}



%

\printaddresses
\end{document}